\documentclass[11pt]{amsart}
\usepackage{amsmath,amssymb,amsopn,amsthm,amscd}

\newtheorem{theorem}{Theorem}[section]
\newtheorem{lemma}[theorem]{Lemma}
\newtheorem{proposition}[theorem]{Proposition}
\newtheorem{corollary}[theorem]{Corollary}

\theoremstyle{definition}
\newtheorem{definition}[theorem]{Definition}
\newtheorem{example}[theorem]{Example}

\theoremstyle{remark}

\begin{document}

\author{Peter Saveliev}
\title{A Lefschetz-type coincidence theorem}

\begin{abstract}
A Lefschetz-type coincidence theorem for two maps $f,g:X\rightarrow Y$ from
an arbitrary topological space to a manifold is given: $I_{fg}=\lambda
_{fg}, $ that is, the coincidence index is equal to the Lefschetz number. It
follows that if $\lambda _{fg}\neq 0$ then there is an $x\in X$ such that $%
f(x)=g(x)$. In particular, the theorem contains well-known coincidence
results for (i) $X,Y$ manifolds, $f$ boundary-preserving, and (ii) $Y$
Euclidean, $f$ with acyclic fibres. It also implies certain fixed point
results for multivalued maps with ``point-like'' (acyclic) and
``sphere-like'' values.
\end{abstract}
\maketitle

\address{Department of Mathematics, University of Illinois at Urbana-Champaign, 1409
West Green Street, Urbana, IL 61801. e-mail: saveliev@member.ams.org}


\section{Introduction.}

A Lefschetz-type coincidence theorem states the following. Given~a pair of
continuous maps $f,g:X\longrightarrow Y,$ the Lefschetz number $\lambda
_{fg} $ of the pair $(f,g)$ is equal to its coincidence index $I_{fg}$,
while $I_{fg}$ is defined in such a way that 
$$
I_{fg}\neq 0\Rightarrow f(x)=g(x)\text{ for some }x\in X. 
$$
Thus, if the Lefschetz number, a computable homotopy invariant of the pair,
is not zero, then there is a coincidence. We now consider two ways to define
the coincidence index in two different settings.

\begin{description}
\item[Case 1]  Let $M_1,M_2$ be closed $n$-manifolds,$\ X$ an open subset
of $M_1,\ N$ an open subset of $M_1,$ $V$ an open subset of $M_2$, $%
f,g:X\longrightarrow V\subset M_2$ maps$,\ \{x\in X:f(x)=g(x)\}\subset
N\subset \overline{N}\subset X\subset M_1.$

Then the {\it coincidence index} $I_{fg}^X$ \cite[p. 177]{Vick} is the image
of the fundamental class $O_{M_1}$ of $M_1$ under the composition:%
$$
\begin{array}{l}
H_n(M_1)^{
\underrightarrow{~inclusion~}}H_n(M_1,M_1\backslash N)^{\underrightarrow{%
~excision~}}H_n(X,X\backslash N) \\ ^{\underrightarrow{(f,g)_{*}}%
}H_n(M_2\times M_2,M_2\times M_2\backslash \delta (M_2))\simeq {\bf Q,}
\end{array}
$$
where $\delta (x)=(x,x).$
\end{description}

This definition represents the original approach to the coincidence problem
for closed manifolds due to Lefschetz \cite{Lef}. It was later generalized
to the case of manifolds with boundary (and a boundary-preserving $f$) by
Nakaoka \cite{Naka}, Davidyan \cite{David0,David}, Mukherjea \cite{Mukh}.

\begin{description}
\item[Case 2]  Let $V$ be an open subset of $n$-dimensional Euclidean
space, $\ f:X\longrightarrow V$ a Vietoris map (i.e., $f^{-1}(y)$ is acyclic
for each $y\in V$), $g:X\longrightarrow K$ a map, $K\subset V$ a finite
polyhedron,$\ N=f^{-1}(K).$

Then the {\it coincidence index} $I(f,g)$ \cite[p. 38]{Gorn} is the image of
the fundamental class $O_K$ of $K$ under the composition:%
$$
H_n(V,V\backslash K)^{\underrightarrow{~f_{*}^{-1}~}}H_n(X,X\backslash N)^{%
\underrightarrow{~(f-g)_{*}~}}H_n({\bf R}^n,{\bf R}^n\backslash \{0\})\simeq 
{\bf Q.} 
$$
\end{description}

This coincidence index has evolved from the Hopf's fixed point index $I_g$ ($%
X=V,f=Id_V$), see Brown \cite[Chapter IV]{Brown}, Dold \cite{Dold0}, 
\cite[VII.5]{Dold}. This approach was developed by Eilenberg and Montgomery 
\cite{EM}, Begle \cite{Begle}, Gorniewicz and Granas \cite{Gorn,GG1,GG2} and
others, see \cite{Gorn} for bibliography. It does not require any knowledge
of the topology of $X$ and, for this reason, is especially well suited for
the study of fixed points of multivalued maps. For an acyclic-valued
multifunction $F:Y\longrightarrow Y,$ we let $X$ be the graph of $F$ and $%
f,g $ be the projections of $X$ on $Y,$ then $f$ is a Vietoris map, and a
coincidence of $(f,g)$ is a fixed point of $F$. This construction can not be
applied to Case 1 because the graph of a multifunction $F:M_2\longrightarrow
M_2$ is not, in general, a manifold.

The restrictions on spaces and maps in Case 1 and Case 2 are necessary to
ensure the existence of an appropriate homomorphism $f_{!}:H(Y)%
\longrightarrow H(X)$, which we shall call here the transfer of $f.$ Then
the Lefschetz number of $\varphi _{fg}=g_{*}f_{!}$ is said to be the
Lefschetz number of the pair $(f,g).$ For Case 1,%
$$
f_{!}=D_1f^{*}D_2^{-1}, 
$$
where $D_1$ and $D_2$ are the Poincare duality isomorphisms for manifolds $%
M_1$ and $M_2.$ For Case 2,%
$$
f_{!}=f_{*}^{-1}, 
$$
with the existence of $f_{*}^{-1}$ guaranteed by the Vietoris-Begle Mapping
Theorem \ref{VietB}.

Until now these two ways to treat the same problem have been studied
separately. In this paper we provide a unified approach. We define the
coincidence index as in Case 1, for arbitrary maps to an $n$-manifold (with
or without boundary), $n\geq 1$, but with no restriction on their domain, as
in Case 2$.$ Roughly, we combine 
$$
\begin{array}{llll}
\text{Case 1:} & n\text{-manifold} & ^{\underrightarrow{~~~\text{any map}%
~~~~~}} & \ n 
\text{-manifold, and} \\ \text{Case 2:} & \text{any space}\  & ^{%
\underrightarrow{~~\text{Vietoris map}~~}} & \text{Euclidean space,} 
\end{array}
$$
into%
$$
\begin{array}{lll}
~\text{any space}~~\  & ^{\underrightarrow{~~~~~\text{any map}~~~~~}} & \ 
\text{manifold.} 
\end{array}
$$
Under the restrictions of Cases 1 and 2, our main theorem reduces to the
results mentioned above (see Sections \ref{Manifolds} and \ref{VietorisMap}%
), but it also applies to the case 
$$
\begin{array}{lll}
~\text{non-manifold}~ & ^{\underrightarrow{~~~~~\text{any map}~~~~~}} & \ 
\text{manifold,} 
\end{array}
$$
as well as%
$$
\begin{array}{lll}
~m\text{-manifold}~ & ^{\underrightarrow{~~~~~\text{any map}~~~~~}} & \ n%
\text{-manifold, }m\neq n 
\end{array}
$$
(see Section \ref{Discussion}). For the sake of simplicity, we limit our
attention to the case when $Y$ is subset of a manifold, although some of the
results can be extended to include spaces as general as ANR's.

An important particular situation when the choice of the transfer is obvious
(see Corollary \ref{maincor1} and Theorem \ref{degree}) occurs if the
conditions below is satisfied:%
$$
{\bf (A)}\ f_{*}:H_n(X,X\backslash N)\rightarrow H_n(V,V\backslash K)\text{
is nonzero} 
$$
(note that, in the case of two $n$-manifolds, the condition simply means
that $f$ has nonzero degree)$.$ This condition can be fairly easily verified
for specific spaces and maps. In particular, we can relax the Vietoris
condition on $f$ or assume that $f$ is a fibration (see Section \ref
{Discussion}). Furthermore, when $g_{*}=0$ (in reduced homology) the
Lefschetz number of $\varphi _{fg}=g_{*}f_{!}$ is equal to $1$, so condition
(A) implies the existence of a coincidence. For example, $f:({\bf D}^2,{\bf S%
}^1)\rightarrow ({\bf S}^2,\{*\})$ has a coincidence with any $g:{\bf D}%
^2\rightarrow {\bf S}^2.$

Our approach seems to be related to a suggestion made by Dold in \cite{Dold2}%
. The subject of his paper is coincidences on ENR$_B$'s, Euclidean
neighborhood retracts over space $B$, and will remain outside the scope of
the present paper. In the end of his paper Dold compares his Theorem 2.1 to
a result that assumes that one of the maps is Vietoris (Case 2): ``It
appears less general than 2.1 because 2.1 makes no acyclicity assumption ...
On the other hand, it has a more general aspect than 2.1 because it does not
assume an actual fibration (or ENR$_B$), only a ``cohomology fibration''
(with ``pointlike'' fibres). This comparison suggests a common
generalization, namely to general {\it cohomology fibrations }...'' (cf.
Section \ref{GenCase2}).

The proofs of our main theorems is self-contained and use some constructions
from Gorniewicz \cite[V.5.1, pp. 38-40]{Gorn} (see also Dold 
\cite[VII.6, pp. 207-211]{Dold}) and Vick \cite[Chapter 6]{Vick}.

The paper is organized as follows. In Section \ref{MainResults} we present
our main results (Theorems \ref{representation} - \ref{main1}) and in
Section \ref{ExTransfer} we prove Theorem \ref{representation}. In Section 
\ref{Manifolds} we prove Theorem \ref{main1} (a Lefschetz-type coincidence
theorem for maps to a manifold with boundary) and obtain Nakaoka's
Coincidence Theorem for boundary-preserving maps between manifolds and
Gorniewicz's Coincidence Theorem for Vietoris maps. Section \ref{Discussion}
is devoted to applications of the main theorem and examples with the
emphasis on situations that are not covered by the two traditions discussed
above. Sections \ref{CoinIndex} - \ref{MainRes} contain the proof of Theorem 
\ref{identity} in a slightly more general setting (for maps to an open
subset of a manifold). In the last two sections we prove a Lefschetz-type
coincidence theorem for generalized Case 2.

\section{Main Results.\label{MainResults}}

Let $E=\{E_q\}$ be a graded ${\bf Q}$-module with%
$$
\dim E_q<\infty ,\ q=0,\ldots ,n,\text{\quad }E_q=0,\ q=n+1,\ldots 
$$
(in other words, $E$ is finitely generated). If $h=\{h_q\}$ is an
endomorphism of $E$ of degree $0$, then the {\it Lefschetz number} $L(h)$ of 
$h$ is defined by%
$$
L(h)=\sum_q(-1)^qtr(h_q), 
$$
where $tr(h_q)$ is the trace of $h_q.$

By $H$ we denote the singular homology and by $\check H$ the \v Cech
homology with compact carriers with coefficients in ${\bf Q}.$\ Throughout
the paper $M$ is an oriented connected compact closed $n$-manifold, $n\geq 0$
(although most results remain valid for a non-orientable $M$ if we take the
coefficient field to be ${\bf Z}_2$)$.$

Let $X$ be a topological space, $N\subset X,$ $M$ be an oriented connected
compact closed $n$-manifold, $(S,\partial S)$ a connected $n$-submanifold
with (possibly empty) boundary $\partial S$ and interior $\stackrel{\circ }{S%
}=S\backslash \partial S$. Let%
$$
f:(X,X\backslash N)\longrightarrow (S,\partial S),\text{\quad }%
g:X\longrightarrow S, 
$$
be continuous maps with $Coin(f,g)=\{x\in X:f(x)=g(x)\}\subset N$. Let%
$$
M^{\times }=(M\times M,M\times M\backslash \delta (M)), 
$$
where $\delta (x)=(x,x)$ is the diagonal map. Then the map $f\times
g:(X,X\backslash N)\times X\longrightarrow M^{\times }$ is well defined.

Fix an element $\mu \in H_n(X,X\backslash N).$

The {\it coincidence index} $I_{fg}$ {\it of the pair} $(f,g)$ (with respect
to $\mu $) is defined by%
$$
I_{fg}=(f\times g)_{*}\delta _{*}(\mu )\in H_n(M^{\times })\simeq {\bf Q}. 
$$
Let $O_S\in H_n(S,\partial S)$ be the fundamental class of $(S,\partial S).$
The {\it transfer of }$f$ (with respect to $\mu $) is the homomorphism $%
f_{!}:H(S)\longrightarrow H(X)$ given by%
$$
f_{!}=(f^{*}D^{-1})\frown \mu , 
$$
where $D:H^{*}(S,\partial S)\rightarrow H(S)$ is the Poincare-Lefschetz
duality isomorphism. Then we define $\lambda _{fg}=L(g_{*}f_{!})$ to be {\it %
the Lefschetz number of the pair }$(f,g)$ (with respect to $\mu $)$.$

The proofs of the three theorems below are located in Sections \ref
{ExTransfer}, \ref{MainRes} and \ref{Manifolds} respectively.

\begin{theorem}
\label{representation}For a pair $f:(X,X\backslash N)\longrightarrow
(S,\partial S),$\quad $g:X\longrightarrow \stackrel{\circ }{S},$%
$$
I_{fg}=I_{*}(Id\otimes g_{*}f_{!})\delta _{*}(O_S), 
$$
where $I:(S,\partial S)\times \stackrel{\circ }{S}\longrightarrow M^{\times }
$ is the inclusion.
\end{theorem}

\begin{theorem}
\label{identity}For any homomorphism $\varphi :H(S)\rightarrow H(\stackrel{%
\circ }{S})$ we have 
$$
L(\varphi i_{*})=I_{*}(Id\otimes \varphi )\delta _{*}(O_S), 
$$
where $i:\stackrel{\circ }{S}\longrightarrow S$ is the inclusion.
\end{theorem}

The following is the main theorem of the paper.

\begin{theorem}[Lefschetz-Type Theorem]
\label{main1}For any pair $f:(X,X\backslash N)\longrightarrow (S,\partial S),
$\ $g:X\longrightarrow S,$ the coincidence index is equal to the Lefschetz
number (with respect to $\mu $):%
$$
I_{fg}=L(g_{*}f_{!}). 
$$
Moreover, if $L(g_{*}f_{!})\neq 0,$ then $(f,g)$ has a coincidence.
\end{theorem}

If $(X,X\backslash N)$ is a manifold with boundary, we get the
Lefschetz-type coincidence theorem for Case 1 by letting $\mu $ be its
fundamental class. To get such a theorem for Case 2 we let $M={\bf R}^n\cup
\{\infty \},\ \mu =f_{*}^{-1}(O_S)$.

Now we consider these results in detail.

\section{The Transfer of a Map and Theorem \ref{representation}. \label
{ExTransfer}}

Recall \cite[p. 156]{Vick} that if $(S,\partial S)$ is a compact oriented $n$%
-manifold, then the Poincare-Lefschetz duality isomorphism%
$$
D:H^{n-k}(S,\partial S)\longrightarrow H_k(S) 
$$
is given by $D(a)=a\frown O_S$. Suppose $f:(S_1,\partial S_1)\longrightarrow
(S_2,\partial S_2),$ where $(S_i,\partial S_i),\ i=1,2,$ are $n$-manifolds,
is a map. Following Vick \cite[Chapter 6]{Vick} we could define $f_{!}$ as
follows. If

$$
D_i:H^{n-k}(S_i,\partial S_i)\longrightarrow H_k(S_i),\text{\quad }i=1,2, 
$$
denote the duality isomorphisms, we let%
$$
f_{!}=D_1f^{*}D_2^{-1}, 
$$
so that $f_{!}$ is the composition of the following maps: 
$$
H_k(S_2)\stackrel{D_2^{-1}}{\longrightarrow }H^{n-k}(S_2,\partial S_2)%
\stackrel{f^{*}}{\longrightarrow }H^{n-k}(S_1,\partial S_1)\stackrel{D_1}{%
\longrightarrow }H_k(S_1). 
$$

Similarly we define $f_{!}$ for $f:(X,X\backslash N)\rightarrow (S,\partial
S)$, where $X$ is an arbitrary topological space: the transfer of $f$ is the
homomorphism $f_{!}:H(S)\longrightarrow H(X)$ given by%
$$
f_{!}=(f^{*}D_2^{-1})\frown \mu , 
$$
where $D_2:H^{*}(S,\partial S)\rightarrow H(S)$ is the Poincare-Lefschetz
duality isomorphism.

To prove Theorem \ref{representation} we will use some arguments from Vick 
\cite[pp. 184-186]{Vick}. Select a basis $\{x_i\}$ for $H^{*}(S)$ and denote 
$\{a_i\}$ the basis for $H(S)$ dual to $\{x_i\}$ under the Kronecker index.
Define a basis $\{x_i^{\prime }\}$ for $H^{*}(S,\partial S)$ by requiring
that $D_2(x_i^{\prime })=a_i$ and let $\{a_i^{\prime }\}$ be the basis for $%
H(S,\partial S)$ dual to $\{x_i^{\prime }\}$ under the Kronecker index. Thus
we have%
$$
\begin{array}{c}
<x_i,a_j>=<x_i^{\prime },a_j^{\prime }>=\delta _{ij}, \\ 
D_2(x_i^{\prime })=x_i^{\prime }\frown O_S=a_i. 
\end{array}
$$
Similarly, select a basis $\{y_i^{\prime }\}$ for $H^{*}(X,X\backslash N)$
and denote by $\{b_i^{\prime }\}$ the basis for $H(X,X\backslash N)$ dual to 
$\{y_i^{\prime }\}$ under the Kronecker index. We define the homomorphism $%
D_1:H^{*}(X,X\backslash N)\rightarrow H(X)$ by $D_1(x)=x\frown \mu $ and we
let $b_i=D_1(y_i^{\prime })$. Next we let $\{y_i\}$ $\subset H^{*}(X)$ be a
collection dual to $\{b_i\}$ under the Kronecker index. Thus we have%
$$
\begin{array}{c}
<y_i,b_j>=<y_i^{\prime },b_j^{\prime }>=\delta _{ij}, \\ 
D_1(\ y_k^{\prime })=y_i^{\prime }\frown \mu =b_i. 
\end{array}
$$

\begin{lemma}
\label{lembas}(cf. Vick \cite[Lemma 6.10, p. 185]{Vick})%
$$
\sum\limits_i(Id\times f_{!})(a_i^{\prime }\times
a_i)=\sum\limits_i(f_{*}\times Id)(b_i^{\prime }\times b_i). 
$$
\end{lemma}

\begin{proof} Since $\{y_i^{\prime }\}$ and $\{a_i^{\prime }\}$ are bases, there
are representations%
$$
\begin{array}{c}
f^{*}(x_i^{\prime })=\sum\limits_k\gamma _{ik}y_k^{\prime }\text{ and }%
f_{*}(b_j^{\prime })=\sum\limits_k\lambda _{kj}a_k^{\prime }\text{.}
\end{array}
$$
Then%
$$
\gamma _{ij}=<\sum_k\gamma _{ik}y_k^{\prime },b_j^{\prime
}>=<f^{*}(x_i^{\prime }),b_j^{\prime }>=<x_i^{\prime },f_{*}(b_j^{\prime
})>=<x_i^{\prime },\sum_k\lambda _{kj}a_k^{\prime }>=\lambda _{ij}, 
$$
so%
$$
f_{*}(b_i^{\prime })=\sum\limits_k\gamma _{ki}a_k^{\prime }. 
$$
Next,%
$$
D_1f^{*}D_2^{-1}(a_i)=D_1f^{*}(x_i^{\prime })=D_1\sum_k\gamma
_{ik}y_k^{\prime }=\sum_k\gamma _{ik}b_k. 
$$
Therefore,%
$$
\begin{array}{ll}
(Id\times f_{!})(a_i^{\prime }\times a_i) & =a_i^{\prime }\times
D_1f^{*}D_2^{-1}(a_i) \\  
& =a_i^{\prime }\times \sum\limits_k\gamma _{ik}b_k \\  
& =\sum\limits_k\gamma _{ik}(a_i^{\prime }\times b_k).
\end{array}
$$
On the other hand,%
$$
\begin{array}{ll}
(f_{*}\times Id)(b_i^{\prime }\times b_i) & =f_{*}(b_i^{\prime })\times b_i
\\  
& =(\sum\limits_k\gamma _{ki}a_k^{\prime })\times b_i \\  
& =\sum\limits_k\gamma _{ki}(a_k^{\prime }\times b_i).
\end{array}
$$
Therefore summation over $i$ produces the same result. \end{proof}

\begin{lemma}
\label{lemdel}(cf. Vick \cite[Lemma 6.11, p. 186]{Vick}) 
$$
\begin{array}{l}
\text{(a) }\delta _{*}(O_S)=\sum\limits_i(a_i^{\prime }\times a_i), \\ 
\text{(b) }\delta _{*}(\mu )=\sum\limits_i(b_i^{\prime }\times b_i).
\end{array}
$$
$\ $
\end{lemma}

\begin{proof} (a) By the K\"unneth formula, $\{a_i^{\prime }\times a_j\}$ is a
basis of $H((S,\partial S)\times S),$ and $\{x_i^{\prime }\times x_j\}$ is
the dual basis of $H^{*}((S,\partial S)\times S)$. Then the identity follows
from the equations:%
$$
\begin{array}{ll}
<x_j^{\prime }\times x_k,\delta _{*}(O_S)> & =<\delta ^{*}(x_j^{\prime
}\times x_k),O_S> \\  
& =<x_j^{\prime }\smile x_k,O_S> \\  
& =<x_k,x_j^{\prime }\frown O_S> \\  
& =<x_k,a_j> \\  
& =\delta _{kj.}
\end{array}
$$

(b) It is clear that $\delta _{*}(\mu )$ belongs to the subspace of $%
H((X,X\backslash N)\times X)$ spanned by $\{b_i^{\prime }\times b_i\}.$ Then
the rest of the proof follows (a). \end{proof}

\begin{proof}[Proof of Theorem \ref{representation}.] We have%
$$
\begin{array}{ll}
I_{fg}=(f\times g)_{*}\delta _{*}(\mu ) & =I_{*}(Id\times g_{*})(f_{*}\times
Id)\sum\limits_i(b_i^{\prime }\times b_i)\ 
\text{by Lemma }\ref{lemdel} \\  & =I_{*}(Id\times g_{*})(Id\times
f_{!})\sum\limits_i(a_i^{\prime }\times a_i)\ 
\text{by Lemma }\ref{lembas} \\  & =I_{*}(Id\times g_{*}f_{!})\delta
_{*}(O_S)\ \text{by Lemma }\ref{lemdel}.
\end{array}
$$
\end{proof}

Thus $I_{fg}$ is the image of $O_S$ under the composition of the following
maps:%
$$
H(S,\partial S)^{\underrightarrow{\ \delta _{*}\ }}H(S,\partial S)\otimes
H(S)^{\underrightarrow{Id\otimes \varphi \ \ }}H(S,\partial S)\otimes H(%
\stackrel{\circ }{S})^{\underrightarrow{\ I_{*}\ }}H(M^{\times }), 
$$
while $\varphi =g_{*}f_{!}$ is defined by the following diagram:%
$$
\begin{array}{ccc}
H^{*}(X,X\backslash N) & ^{\ \underleftarrow{\quad \ \ \ \ \ f^{*}\ \quad }%
}\  & H^{*}(S,\partial S) \\ 
\downarrow ^{\frown \mu } &  & \downarrow ^{D_2} \\ 
H(X) & ^{\underrightarrow{\ \ \ \ \ \ g_{*}\ \ \ \ \ \ }} & H(S). 
\end{array}
$$

It is worth mentioning that Nakaoka \cite{Naka} defines the {\it coincidence
transfer homomorphism} $\tau _{fg}$ for a pair of fibre-preserving maps $%
f,g:E\rightarrow E^{\prime }$, where $E,E^{\prime }$ are manifolds with
boundary. If the base is trivial then, according to \cite[Theorem 5.1 (ii)]
{Naka}, $\tau _{fg}$ is related to $\varphi $ as follows:%
$$
\tau _{fg}(1)=L(\varphi _{fg}). 
$$

\section{The Main Theorem and Cases 1 and 2. \label{Manifolds}}

The sum of Theorems \ref{representation} and \ref{identity} gives us a
Lefschetz-type coincidence theorem for $g:X\rightarrow \stackrel{\circ }{S}$%
. To prove the statement for $g:X\rightarrow S$ (Theorem \ref{main1}) we
need the following fact.

\begin{proposition}[Topological Collaring Theorem]
\cite[Theorem 5.2, p. 154]{Vick}Let $(S,\partial S)$ be a manifold. Then
there is a manifold $(T,\partial T)$ obtained from $(S,\partial S)$ by
attaching a ``collar'':%
$$
T=S\cup (\partial S\times [0,1]). 
$$
\end{proposition}

Next we restate and prove Theorem \ref{main1}.

\begin{theorem}
\label{LefManifold}For any pair $f:(X,X\backslash N)\longrightarrow
(S,\partial S),$\ $g:X\longrightarrow S,$ with $Coin(f,g)\subset N$, the
coincidence index is equal to the Lefschetz number (with respect to $\mu $):%
$$
I_{fg}=L(g_{*}f_{!}). 
$$
Moreover, if $L(g_{*}f_{!})\neq 0,$ then $(f,g)$ has a coincidence.
\end{theorem}

\begin{proof} We attach ``collars'' to $X$ and $S$ as follows. Let%
$$
Z=X\cup ((X\backslash N)\times [0,1]), 
$$
such that $X\cap ((X\backslash N)\times [0,1])=(X\backslash N)\times \{0\}$.
And we assume that according to the proposition above 
$$
\ T=S\cup (\partial S\times [0,1])\subset M 
$$
is an $n$-manifold. Then we define $G:Z\rightarrow \stackrel{\circ }{T}$ by
(cf. \cite{David}) 
$$
G=jgr, 
$$
where $j:S\longrightarrow \stackrel{\circ }{T}$ is the inclusion, $%
r:Z\longrightarrow X$ is the retraction. We also define $F:(Z,(X\backslash
N)\times \{1\})\longrightarrow (T,\partial T)$ by%
$$
\begin{array}[t]{ll}
F(x,t)=(f(x),t) & \text{\ if }(x,t)\in (X\backslash N)\times [0,1], \\ 
F(x)=f(x) & \text{\ if }x\in X.
\end{array}
$$
From Theorems \ref{representation} and \ref{identity} we have%
$$
I_{FG}=L(G_{*}F_{!}i_{*}), 
$$
where $i:\ \stackrel{\circ }{T}\rightarrow T$ is the inclusion. The
inclusions and retractions induce isomorphisms, $F_{*}=f_{*}$, $G_{*}=g_{*}$%
, so $L(G_{*}F_{!}i_{*})=L(g_{*}f_{!}).$ Next, the following diagram
commutes:%
$$
\begin{array}{ccc}
(X,X\backslash N) & ^{\underrightarrow{\quad (f,g)\quad }} & M^{\times } \\ 
~~~~~~\uparrow ^r\quad  &  & ~~~~\parallel \ \ \ \  \\ 
(Z,(X\backslash N)\times \{1\}) & ^{\underrightarrow{\quad (F,G)\quad }} & 
M^{\times },
\end{array}
$$
which means that $I_{FG}=I_{fg}.$ Thus,%
$$
I_{fg}=I_{FG}=L(G_{*}F_{!}i_{*})=L(g_{*}f_{!}). 
$$
\end{proof}

Halpern \cite{Halpern} proves a Lefschetz-type coincidence theorem in an
even more general situation: he considers $f,g:X\rightarrow Y,$ where both $%
X $ and $Y$ are arbitrary topological spaces. His Lefschetz number is the
Lefschetz number of the homomorphism $\varphi :Y\rightarrow Y$ given by%
$$
\varphi (z)=g_{*}((f^{*}(b/z)\frown a) 
$$
for some elements $a\in H_n(X)$ and $b\in H^n(Y\times Y),$ and proves that $%
L(\varphi )\neq 0$ implies that $Coin(f,g)\neq \emptyset .$ To compare his
result with ours, observe first that he does not define the coincidence
index, which has an independent interest, and second his theorem does not
include the Brouwer fixed point theorem.

\subsection{Case 1.}

To get a Lefschetz-type coincidence theorem for Case 1, i.e., when $%
(X,X\backslash N)$ is an $n$-manifold, we simply let $\mu $ be its
fundamental class.

\begin{corollary}
\label{Lefmani}Let $(S_1,\partial S_1),(S_2,\partial S_2)$ be oriented
compact connected $n$-manifolds, and let%
$$
f:(S_1,\partial S_1)\longrightarrow (S_2,\partial S_2),\text{\quad }%
g:S_1\longrightarrow S_2 
$$
be continuous maps$.$ If $Coin(f,g)\cap \partial S_1=\emptyset ,$ then the
coincidence index with respect to $\mu =O_{S_1}$ is equal to the Lefschetz
number:%
$$
I_{fg}=L(g_{*}f_{!}) 
$$
(here $f_{!}$ is defined via Poincare duality for $S_1$ and $S_2$).
Moreover, if $L(g_{*}f_{!})\neq 0$ then $(f,g)$ has a coincidence.
\end{corollary}

Several authors have dealt with a Lefschetz-type coincidence theorem for
manifolds with boundary. Corollary \ref{Lefmani} provides little additional
information in comparison to these results but still has certain advantages.

The Lefschetz-Nakaoka Coincidence Theorem \cite[Theorem 3.2]{BS} (it is
Lemma 8.1 combined with Theorem 5.1 of Nakaoka \cite{Naka}) is identical to
our Theorem \ref{Lefmani} but applies only to manifolds with nonempty
boundary. The reason for this is that in \cite{Naka} manifolds with boundary
are ``doubled'' (two copies are glued together along the boundary) and then
Nakaoka's Lefschetz-type coincidence theorem for closed manifolds is
applied. Of course, the case of empty boundary follows from the classical
Lefschetz coincidence theorem \cite[Theorem 6.13]{Vick}, but the case $%
\partial S_1=\emptyset ,$ $\partial S_2\neq \emptyset $ is still excluded.
Bredon \cite[VI.14]{Bredon} also considers manifolds with empty and nonempty
boundary separately. Theorem of Davidyan \cite{David} and Theorem 2.1 of
Mukherjea \cite{Mukh} use collaring instead of doubling, so they can be
specialized to manifolds with empty boundary. But they do not prove that the
coincidence index is equal to the Lefschetz number (Davidyan \cite{David0}
proves this identity only for manifolds without boundary).

Therefore Corollary \ref{Lefmani} is of some interest, because it opens a
possibility of defining a coincidence index for all manifolds with boundary,
empty or not. Such an index may be used for a unified Nielsen coincidence
theory, see Brown and Schirmer \cite{BS,BS1}, where boundaries are required
to be nonempty.

\subsection{Case 2.}

We can also obtain a Lefschetz-type coincidence theorem for Case 2 (cf.
Theorem \ref{GornCoin}) with an additional assumption.

Recall \cite[p. 13]{Gorn} that a map $f:(X,X_0)\rightarrow (Y,Y_0)$ is said
to be {\it Vietoris }if (i) $f$ is proper, i.e., $f^{-1}(B)$ is compact for
any compact $B\subset Y,$ (ii) $f^{-1}(Y_0)=X_0,$ (iii) the set $f^{-1}(y)$
is acyclic with respect to the \v Cech homology for every $y\in Y.$

\begin{proposition}[Vietoris-Begle Theorem]
\label{VietB}\cite[p. 14]{Gorn} If $f:(X,X_0)\rightarrow (Y,Y_0)$ is a
Vietoris map, then $f_{*}:\check H(X,X_0)\rightarrow \check H(Y,Y_0)$ is an
isomorphism.
\end{proposition}

\begin{corollary}
\label{GornCoin1}Suppose $X$ is a topological space, $S\subset {\bf R}^n$ is
a compact $n$-manifold. Suppose%
$$
f,g:X\rightarrow \stackrel{\circ }{S} 
$$
are two continuous maps such that $f$ is Vietoris. Then if $%
L(g_{*}f_{*}^{-1})\neq 0$ with respect to the \v Cech homology over ${\bf Q}$%
, then the pair $(f,g)$ has a coincidence.
\end{corollary}

\begin{proof} We put $\ \mu =f_{*}^{-1}(O_S)$ and apply Theorem \ref{LefManifold}.
Then $I_{fg}=L(g_{*}f_{!}).$ But by Theorem \ref{VietB}, $f_{*}^{-1}$
exists, hence by the proposition below, $f_{!}=f_{*}^{-1}$. Therefore $%
I_{fg}=L(g_{*}f_{*}^{-1}).$  \end{proof}

\begin{proposition}
\label{f-inverse}(cf. \cite[Proposition VI.14.1 (6), p. 394]{Bredon}) If $%
f_{*}(\mu )=O_S$ then $f_{*}f_{!}=Id.$
\end{proposition}

\section{Corollaries and Examples.\label{Discussion}}

Suppose $(X,X^{\prime })$ is a topological space, $(S,\partial S)$ is an
oriented compact connected $n$-manifold. Consider the following condition:

\begin{description}
\item[(A)]  $f_{*}:H_n(X,X^{\prime })\rightarrow H_n(S,\partial S)$ is a
nonzero homomorphism.
\end{description}

First we will consider analogues of some well-known theorems about maps
between manifolds without the assumption that the domain of the maps is a
manifold. The following is a generalization of Theorem 2.2 of Mukherjea \cite
{Mukh}.

\begin{corollary}
\label{maincor1}Suppose $f:(X,X^{\prime })\longrightarrow (S,\partial S)$
satisfies condition (A) and $g:X\longrightarrow S$ induces $g_{*}=0$ (in
reduced homology). Then $(f,g)$ has a coincidence.\ 
\end{corollary}

\begin{proof} First we select $\mu \in H_n(X,X^{\prime })$ such that $f_{*}(\mu
)=O_S$. Then by Proposition \ref{f-inverse}, $f_{*}f_{!}=Id$. Then $%
g_{*}f_{!}:H_i(S,\partial S)\rightarrow H_i(S,\partial S)$ is non-zero for $%
i=0$ and zero for $i\neq 0$. Therefore $L(g_{*}f_{!})=1$, so there is a
coincidence by Theorem \ref{main1}. \end{proof}

Corollary \ref{maincor1} allows us to use a version of the Vietoris Theorem
stronger than Proposition \ref{VietB} (see also Section \ref{Discussion}):
even if $H_n(f^{-1}(x))\neq 0$ for some $x$, while $H_i(f^{-1}(x))=0$ for
all $x$ and $0\leq i<n$, we still have an epimorphism $f_{*}:H_n(X,X%
\backslash N)\rightarrow H_n(S,\partial S).$ Therefore condition (A) is
satisfied. For other versions of the Vietoris Theorem see \cite{Krys}.

The following is a generalization of the Kronecker theorem: a map with
nonzero degree is onto.

\begin{corollary}
If $f:(X,X^{\prime })\rightarrow (S,\partial S)$ satisfies condition (A)
then $f$ is onto.
\end{corollary}

\begin{proof} For a given $y\in S,$ we define $g:X\rightarrow S$ by $g(x)=y,$ for
all $x\in X$. Therefore $g_{*}:H(X)\rightarrow H(S)$ is a zero homomorphism,
hence there is a coincidence by the previous corollary. Therefore $y\in $ $%
f(X)$, so $f(X)=S$. \end{proof}

A map $f:(S_1,\partial S_1)\longrightarrow (S_2,\partial S_2)$, where $%
(S_i,\partial S_i),\ i=1,2,$ are manifolds, is called {\it %
coincidence-producing} \cite[Section 7]{BS} if every map $%
g:S_1\longrightarrow S_2$ has a coincidence with $f$. Brown and Schirmer 
\cite[Theorem 7.1]{BS} showed that if $S_2$ is acyclic, $n\geq 2,$ then $f$
is coincidence-producing if and only if $f_{*}:H_n(S_1,\partial
S_1)\rightarrow H_n(S_2,\partial S_2)$ is nonzero. We call a map $%
f:(X,X^{\prime })\rightarrow (S,\partial S)$ {\it weakly
coincidence-producing} if every map $g:X\rightarrow S$ with $g_{*}=0$ has a
coincidence with $f.$ Then Corollary \ref{maincor1} takes the following form.

\begin{corollary}
\label{nonV}If $f:(X,X^{\prime })\rightarrow (S,\partial S)$ satisfies
condition (A) then $f$ is weakly coincidence-producing.
\end{corollary}

For an acyclic $S$, this is a generalization of the ``if'' part of the
Brown-Schirmer statement.

We conclude with a few examples of applications of Corollary \ref{nonV}.
These examples are not included in either Case 1 or Case 2.

\subsection{Manifolds.}

It is hard to come by an example of coincidences that does not involve
manifolds. Yet we can consider a pair $(X,X^{\prime })$ such that $X$ is a
manifold (possibly with boundary) but $(X,X^{\prime })$ is not a manifold
with boundary, i.e., $X^{\prime }$ is not the boundary of $X$ (nor
homotopically equivalent to it). This provides a setting not included in
Case 1.

\begin{example}
Let $f:({\bf D}^2,\partial e\cup \partial e^{\prime })\rightarrow ({\bf D}^2,%
{\bf S}^1),$ where $e$ and $e^{\prime }$ are disjoint cells in ${\bf D}^2$,
be a map.
\end{example}

Then it is a matter of simple computation to check whether condition (A) is
satisfied. For examples of acyclic manifolds, see Brown and Schirmer 
\cite[Section 7]{BS}.

\begin{example}
Let $f:({\bf I},\partial {\bf I})\times {\bf S}^1\rightarrow ({\bf D}^2,{\bf %
S}^1\cup \{0\}),\ {\bf I}=[0,1]$, be the map that takes by identification $%
\{0\}\times {\bf S}^1$ to $\{0\}.$
\end{example}

It is clear that condition (A) is satisfied for $X^{\prime }=\{1\}\times 
{\bf S}^1$, therefore by Corollary \ref{nonV}, any map homotopic to $f$ has
a coincidence with a map $g$ if $g_{*}=0$. But $f^{-1}(0)={\bf S}^1$ is not
acyclic, so this example is not covered by Case 2 and Gorniewicz's Theorem.
On the other hand, even though $f$ maps a manifold to a manifold, it does
not map boundary to boundary. Therefore, Case 1 and the results discussed in
Section \ref{Manifolds} do not include this example.

In a similar fashion we can show that the projection of the torus ${\bf T}^2$
on the circle ${\bf S}^1$ is a weakly coincidence-producing map. This is an
example of a map between manifolds of different dimensions. For a negative
example of this kind, take the Hopf map $f:{\bf S}^3\rightarrow {\bf S}^2$,
then for any $g,$ the Lefschetz number and the coincidence index of the pair 
$(f,g)$ are equal to zero.

\subsection{Non-manifolds.}

Let $E$ be a space that is not acyclic and not a manifold, e.g., the
``figure eight''. Consider the projection $f:(X,X^{\prime })=({\bf I}%
,\partial {\bf I})\times E\rightarrow ({\bf I},\partial {\bf I}),{\bf I}%
=[0,1],$ onto the first coordinate. Then $f$ clearly satisfies (A). Thus $f$
is weakly coincidence-producing by Corollary \ref{nonV}. Observe that $%
f^{-1}(x)=E$ is not acyclic and $X$ is not a manifold.

A relevant example is given in Kahn \cite{Kahn}. He constructed an infinite
dimensional acyclic space $X$ and an essential map $f:X\rightarrow {\bf S}%
^3. $ Then $f$ satisfies condition (A) and, since any $g:X\rightarrow {\bf S}%
^3$ induces a zero homomorphism in reduced homology, it follows that $f$ is
weakly coincidence-producing by Corollary \ref{nonV}.

\subsection{Fibrations and $UV^n$-maps.}

\begin{corollary}
Suppose $X$ is a topological space, $M$ is an oriented compact closed $(n-1)$%
-connected $n$-manifold, $f:X\rightarrow M$ is a map, and

\begin{description}
\item[(A$^{\prime }$)] $f_{\#}:\pi _n(X)\rightarrow \pi _n(M)$ is onto$.$
\end{description}

\noindent Then $f$ is weakly coincidence producing.
\end{corollary}

\begin{proof} As $M$ is $(n-1)$-connected, the Hurewicz homomorphism $h_n:\pi
_n(M)\rightarrow H_n(M)$ is onto \cite[p. 488]{Bredon}. Hence $f_{*}$ is
onto and $(A)$ is satisfied. \end{proof}

Condition (A$^{\prime }$) holds when $f$ is a fibration with $\pi
_{n-1}(f^{-1}(y))=0$ (it follows from the homotopy sequence of the fibration 
\cite[Theorem VII.6.7, p. 453]{Bredon}).

Condition (A$^{\prime }$) also holds when $f$ is onto and for each $y\in Y,\
f^{-1}(y)$ has the $UV^n${\it -property} for each $n$ (see \cite[Section 4]
{Krys} and its bibliography): for any neighborhood $U$ of $f^{-1}(x)$ there
is a neighborhood $V\subset U$ such that any singular $k$-sphere in $V$ is
inessential, $0\leq k\leq n$. Then Theorem \ref{main} gives a version of
Theorem 1.2 of Gutev \cite{Gutev}. For related results see also 
\cite[Section 2]{Dran}.

\subsection{$m$-Acyclic Maps.}

A multivalued map $\Phi :Y\rightarrow Y$ is called $m${\it -acyclic}, $m\geq
1$, if for each $x\in Y,$ $\Phi (x)$ consists of exactly $m$ acyclic
components. Schirmer \cite{Schirmer} proved that if $Y$ is locally connected
and simply connected then the graph $X$ of $\Phi $ is a disjoint union of
graphs $X_i$ of $m$ acyclic maps $\Phi _i,\ i=1,2,...,m$. Then the
projections $f_i:X_i\rightarrow Y$ on the first coordinate induce
isomorphisms, therefore condition (A) holds for $f$ the projection of $X$
onto $Y$.

Patnaik \cite{Patnaik} defines an $m${\it -map} $\Phi $ as a multifunction
such that $\Phi (x)$ contains exactly $m$ points for each $x$. Then he
considers the Lefschetz number of a homomorphism similar to $g_{*}f_{!}$,
although it is not clear how it is related to ours.

\subsection{Spherical Maps.}

Continuing the discussion in the beginning of this section, what if the
acyclicity condition for $f$ fails at $(n-1)$ degree? Then there is no
version of the Vietoris Theorem available to ensure condition (A).

\begin{definition}
(cf. \cite{ON,Gorn1,Dawid}) Let $B(A),A\subset {\bf R}^n,$ denote the
bounded component of ${\bf R}^n\backslash A.$ A closed-valued u.s.c. map $%
\Phi :{\bf D}^n\rightarrow {\bf D}^n$ is called $(n-1)${\it -spherical}, $n>1
$, if

\begin{description}
\item[(i)] for every $x\in {\bf D}^n$, $H(\Phi (x))=H({\bf S}^{n-1})$ or $%
H(point),$

\item[(ii)] for every $x\in {\bf D}^n$, if $x\in B(\Phi (x))$ then there
exists an $\varepsilon $-neighborhood $O_\varepsilon (x)$ of $x$ such that $%
x^{\prime }\in B(\Phi (x^{\prime }))$ for each $x^{\prime }\in O_\varepsilon
(x).$
\end{description}
\end{definition}

\begin{corollary}
An $(n-1)$-spherical map $\Phi :{\bf D}^n\rightarrow {\bf D}^n,\ n>1,$ has a
fixed point.
\end{corollary}

\begin{proof} We notice, first, that if $\Phi $ has no fixed points and there are
no points $x$ such that $x\in B(\Phi (x)),$ then by replacing $\Phi (x)$
with $\Phi ^{\prime }(x)=\Phi (x)\cup B(\Phi (x))$ we obtain an acyclic
multifunction without fixed points. Therefore we suppose that such an $x$
exists and for simplicity assume that it is $0$. Now, if $0$ is not a fixed
point, then from the upper semicontinuity of $\Phi $ and (ii) above, it
follows that there is an $\varepsilon >0$ such that 
$$
|x|<\varepsilon \Rightarrow x\in B(\Phi (x))\text{ and }|\Phi
(x)|>2\varepsilon . 
$$
Let $X$ be the graph of $\Phi $, $K=\{x:|x|\geq 2\varepsilon \},\ f,g$
projections of $X,\ X^{\prime }=f^{-1}({\bf D}^n\backslash K)$. One can see
that $f$ is essentially the same as the projection of $({\bf D}^n,{\bf S}%
^{n-1})\times {\bf S}^1$ onto $({\bf D}^n,{\bf S}^{n-1})$, and, therefore,
induces a surjection 
$$
f_{*}:H_n(({\bf D}^n,{\bf S}^{n-1})\times {\bf S}^1)\rightarrow H_n({\bf D}%
^n,{\bf S}^{n-1}). 
$$
Hence condition (A) is satisfied, so by Corollary \ref{nonV}, $\Phi $ has a
fixed point. \end{proof}

\begin{example}
(O'Neill \cite{ON}) Let $\Phi :{\bf D}^2\rightarrow {\bf D}^2$ be given by%
$$
\Phi (x)=\{y\in {\bf D}^2:|y-x|=\rho (x)\}\cup \{y\in {\bf S}^1:|y-x|>\rho
(x)\}, 
$$
where $\rho (x)=1-|x|+|x|^2,\ x\in {\bf D}^2.$
\end{example}

This example shows that condition (ii) of the above definition is necessary
for existence of a fixed point.

\section{The Coincidence Index of a Pair of Maps.\label{CoinIndex}}

The next three sections are devoted to the proof of Theorem \ref{identity},
which will be carried out in a setting slightly more general than that of
Section \ref{MainResults}.

Consider the sets $K\subset V\subset M.$ Assume $(M,V,M\backslash K)$ is an
excisive triad, i.e., the inclusion $j:(V,V\backslash K)\longrightarrow
(M,M\backslash K)$ induces an isomorphism in homology. Let $%
i:K\longrightarrow V,\ I:M\times K\longrightarrow M^{\times }$ be the
inclusions$.$

We start by recalling some facts about manifolds. The following definition
and propositions are taken from Vick \cite[Chapter 5]{Vick}.

\begin{proposition}
\label{manifolds}\cite[Corollary 5.7, p. 136]{Vick} For each $p\in M,$ the
homomorphism%
$$
i_{p*}:H_n(M)\longrightarrow H_n(M,M\backslash \{p\})=T_p\simeq {\bf Q} 
$$
where $i_p:M\longrightarrow (M,M\backslash \{p\})$ is the inclusion$,$ is an
isomorphism.
\end{proposition}

\begin{definition}
\cite[p. 139]{Vick} {\it The fundamental class} {\it of }$M$ is an element $%
z\in H_n(M)$ such that%
$$
i_{p*}:H_n(M)\longrightarrow H_n(M,M\backslash \{p\})=T_p 
$$
has $i_{p*}(z)$ a generator of $T_p$ for each $p\in M.$
\end{definition}

\begin{proposition}
\label{lpsp}\cite[Lemma 5.12, p. 143]{Vick} 
$$
\xi :{\bf Q}\simeq H_0(M)\simeq H_n(M^{\times }) 
$$
under the homomorphism $\xi $ sending the $0$-chain represented by $p\in M$
into the relative class represented by $l_{p*}(s(p)),$ where%
$$
l_{p*}:H_n(M,M\backslash \{p\})\longrightarrow H_n(M^{\times }) 
$$
is induced by $l_p(x)=(x,p),\ x\in M,$ and $s:M\longrightarrow T=\cup _{p\in
M}T_p$ is the orientation map of $M$, with $s(p)$ a generator of $T_p$, for
each $p\in M.$
\end{proposition}

\begin{proposition}
\label{ip(z)-sp}\cite[Theorem 5.10, p. 140]{Vick} If $s$ is the orientation
map then there is a unique fundamental class $O_M\in H_n(M)$ such that $%
i_{p*}(O_M)=s(p)$ for each $p\in M.$
\end{proposition}

Consider%
$$
M\stackrel{k}{\longrightarrow }(M,M\backslash K)\stackrel{j}{\longleftarrow }%
(V,V\backslash K), 
$$
where $k$ is the inclusion.

\begin{definition}
\label{def-fund}(cf. \cite[p. 16]{Gorn}, \cite[p. 192]{Dold}) {\it The} {\it %
fundamental class} $O_K$ {\it of the pair} $(V,V\backslash K)$ is defined by%
$$
O_K=j_{*n}^{-1}k_{*n}(O_M)\in H_n(V,V\backslash K). 
$$
\end{definition}

Let $X$ be a topological space, $N$ a subset of $X,$ and let%
$$
f,g:X\longrightarrow V, 
$$
be continuous maps. Suppose%
$$
Coin(f,g)\subset N, 
$$
then the map $f\times g:(X,X\backslash N)\times X\longrightarrow M^{\times }$
is well defined.

Fix an element $\mu \in H_n(X,X\backslash N).$

\begin{definition}
\label{defcoin} The {\it coincidence index} $I_{fg}$ {\it of the pair} $(f,g)
$ (with respect to $\mu $) is defined by%
$$
I_{fg}=(f\times g)_{*}\delta _{*}(\mu )\in H_n(M^{\times })\simeq {\bf Q}. 
$$
\end{definition}

In the setting of Case 1 this definition turns into the usual one (cf. 
\cite[p. 177]{Vick}). Another observation (due to the referee): as $I_{fg}$
is defined via a homomorphism from $H_n(X,X\backslash N)$ to ${\bf Q,}$ it
is an element of $H^n(X,X\backslash N)$.

This definition also includes the coincidence index for Case 2, as given in
the proposition below.

\begin{proposition}
\label{I(fg)}Suppose $U$ is a open subset of $n$-dimensional Euclidean
space, $\ f:X\longrightarrow U$ a Vietoris map, $g:X\longrightarrow K$ a
map, where $K\subset U$ is a finite polyhedron. Let $N=f^{-1}(K),\ \mu
=f_{*}^{-1}(O_K)\in H_n(X,X\backslash N).$ Then $I(f,g)=I_{fg}.$
\end{proposition}

\begin{proof} We identify ${\bf R}^n$ with a hemisphere of $M={\bf S}^n$. Then from
Proposition \ref{lpsp} it follows that $H_0({\bf R}^n)\simeq H_n({\bf R}%
^n\times {\bf R}^n,{\bf R}^n\times {\bf R}^n\backslash \delta ({\bf R}^n))$
under the homomorphism $\xi $ sending the $0$-chain represented by $p\in 
{\bf R}^n$ into the relative class represented by $l_{p*}(s(p))$. Therefore
we have a commutative diagram%
$$
\begin{array}{ccccc}
H_0({\bf R}^n) & ^{\underrightarrow{\qquad \xi \qquad }} & H_n(({\bf R}%
^n)^{\times }) & ^{\underrightarrow{\qquad d_{*}\qquad }} & H_n(
{\bf R}^n,{\bf R}^n\backslash \{0\}) \\ \downarrow  &  & \downarrow  & 
\nwarrow ^{(f,g)_{*}} & \quad \uparrow ^{(f-g)_{*}} \\ 
H_0(M) & ^{\underrightarrow{\qquad \xi \qquad }} & H_n(M^{\times }) & ^{%
\underleftarrow{\qquad (f,g)_{*}\qquad }} & H_n(X,X\backslash N),
\end{array}
$$
where the two vertical arrows are isomorphisms induced by the inclusions, $%
d(x,y)=x-y$. Since $d_{*}$ is an isomorphism \cite[Lemma VII.4.13, p. 200]
{Dold}, it follows that 
$$
I_{fg}=(f,g)_{*}(\mu )=(f-g)_{*}(\mu )=(f-g)_{*}f_{*}^{-1}(O_K)=I(f,g). 
$$
\end{proof}

To justify its name, a coincidence index should satisfy the property below.

\begin{lemma}
\label{indcoin} If $I_{fg}\neq 0$ (with respect to some $\mu $) then the
pair $(f,g)$ has a coincidence.
\end{lemma}

\begin{proof} Suppose not, then $C=Coin(f,g)=\emptyset .$ Hence $H_n((X,X\backslash
C)\times X)=0.$ But the following diagram is commutative:%
$$
\begin{array}{ccc}
H_n((X,X\backslash N)\times X) & ^{\underrightarrow{\qquad (f\times
g)_{*}\qquad }} & {\bf Q} \\ \downarrow ^{k_{*}} &  & \Vert  \\ 
H_n((X,X\backslash C)\times X) & ^{\underrightarrow{\qquad (f\times
g)_{*}\qquad }} & ~{\bf Q,}
\end{array}
$$
where $k$ is the inclusion, so $I_{fg}=0$. \end{proof}

\section{Generalized Dold's Lemma.\label{Dold}}

In this section we obtain a generalization of Dold's Lemma 
\cite[Lemma VII.6.13, p. 210]{Dold} (see also \cite[p. 153]{Brown}), which
is necessary for our definition of the coincidence index. It is proved for
singular homology, but when $V$ is open and $K$ is compact, we can interpret
this result for \v Cech homology with compact carriers, as in 
\cite[I.5.6, p. 17]{Gorn}.

We define the following functions:

{\it the transposition} $t:V\times K\longrightarrow K\times V$ by%
$$
t(x,y)=(y,x); 
$$

{\it the scalar multiplication} $m:{\bf Q\otimes }H(V)\longrightarrow H(V)$
by%
$$
m(r\otimes v)=r\cdot v; 
$$

{\it the tensor multiplication }$O_K^{\times }:H(K)\longrightarrow
H(V,V\backslash K)\otimes H(K),O_M^{\times }:H(K)\longrightarrow H(M)\otimes
H(K)$ by

$$
O_K^{\times }(v)=O_K\otimes v,\text{\quad }O_M^{\times }(v)=O_M\otimes v; 
$$

{\it the projection} $P:H(M^{\times })\longrightarrow H_n(M^{\times })$ by%
$$
P(q)=q_n,\text{\quad if }q=\sum_kq_k,\ q_k\in H_k(M^{\times }). 
$$

\begin{lemma}
\label{prelD}(cf. Dold \cite[Lemma VII.6.14, p. 210]{Dold}) Suppose that $V$
is an ANR. Then the maps%
$$
\psi _0,\psi _1:(V,V\backslash K)\times K\longrightarrow M^{\times }\times V%
\text{ given by} 
$$
$$
\begin{array}{l}
\psi _0(v,k)=(v,k,v), \\ 
\psi _1(v,k)=(v,k,k),\qquad v\in V,\ k\in K,
\end{array}
$$
induce the same homomorphism in homology:%
$$
\psi _{0*}=\psi _{1*}. 
$$
\end{lemma}

\begin{proof} Let $Q={\bf I}\times D\cup \{0,1\}\times V\times K\subset {\bf I}%
\times V\times K$, where ${\bf I}=[0,1],\ D=\{(v,k)\in V\times K:v=k\}$ is
the diagonal of $V\times K.$ Note that $D$ is closed since $V$ is Hausdorff.
Therefore $Q$ is also closed. Consider a function $\alpha :Q\longrightarrow V
$ given by%
$$
\begin{array}{l}
\alpha (0,v,k)=v, \\ 
\alpha (1,v,k)=k, \\ 
\alpha (t,k,k)=k,\qquad \qquad \qquad v\in V,k\in K,t\in {\bf I}.
\end{array}
$$
Clearly $\alpha $ is continuous. Then, since $Q$ is a closed subset of ${\bf %
I}\times V\times K$ and $V$ is an ANR, there is an extension of $\alpha $ to
a neighborhood of $Q$. And since $Q$ contains ${\bf I}\times D$, we assume
that $\alpha $ is now defined on ${\bf I}\times W$, where $W$ is an open
neighborhood of $D$ in $V\times K.$ Suppose maps%
$$
\begin{array}{l}
\eta _i:(W,W\backslash D)\longrightarrow M^{\times }\times V, \\ 
\varphi _i:(V\times K,(V\times K)\backslash D)\longrightarrow M^{\times
}\times V,\qquad i=0,1,
\end{array}
$$
are given by the same formulas as $\psi _i:$%
$$
\begin{array}{l}
\eta _0(v,k)=\varphi _0(v,k)=(v,k,v), \\ 
\eta _1(v,k)=\varphi _1(v,k)=(v,k,k).
\end{array}
$$
Then $\eta _0$ and $\eta _1$ are homotopic:%
$$
\eta _t(v,k)=(v,k,\alpha (t,v,k)),\text{\quad }0\leq t\leq 1. 
$$
Consider the following commutative diagram for $i=0,1$:%
$$
\begin{array}{rrr}
(V,V\backslash K)\times K &  &  \\ 
\downarrow ^j & \stackrel{\psi _i}{\searrow } &  \\ 
(V\times K,(V\times K)\backslash D) & \stackrel{\varphi _i}{\longrightarrow }
& M^{\times }\times V. \\ 
\uparrow ^{j^{\prime }} & \stackrel{\eta _i}{\nearrow } &  \\ 
(W,W\backslash D) &  & 
\end{array}
$$
where $j,j^{\prime }$ are inclusions$.$ Since $W$ is open and $D$ is closed, 
$j_{*}^{\prime }$ is an isomorphism by excision. We also know that $\eta
_0\sim \eta _1,$ so $\eta _{0*}=\eta _{1*}.$ Therefore $\varphi
_{0*}=\varphi _{1*}.$ And since $\psi _i=\varphi _ij$, we finally conclude
that $\psi _{0*}=\psi _{1*}.$ \end{proof}

\begin{theorem}[Generalized Dold's Lemma]
\label{MDold}(cf. \cite[Lemma VII.6.13, p. 210]{Dold}) Suppose that $K$ is
arcwise connected and the map $\Phi :H(K)\longrightarrow H(V)$ is given as a
composition of the following homomorphisms:%
$$
\begin{array}{l}
\Phi :H_i(K)^{
\underrightarrow{~O_K^{\times }~~}}H_{n+i}((V,V\backslash K)\times K)^{%
\underrightarrow{~(\delta \times Id)_{*}~~}}H_{n+i}((V,V\backslash K)\times
V\times K) \\ ^{
\underrightarrow{~(Id\times t)_{*}~~}}H_{n+i}((V,V\backslash K)\times
K\times V)^{\underrightarrow{~(I\times Id)_{*}~~}}H_{n+i}(M^{\times }\times
V) \\ ^{\underrightarrow{~P\otimes Id~~}}H_n(M^{\times })\otimes H_i(V)%
\stackrel{m}{\longrightarrow }H_i(V).
\end{array}
$$
Then%
$$
\Phi =i_{*}. 
$$
\end{theorem}

\begin{proof} Consider the following diagram:%
$$
\begin{array}{ccccccc}
&  & H(M\times K) & ^{\underrightarrow{\ \quad \psi _{*}\quad \ }} & 
H(M\times K\times V) &  &  \\  
& \stackrel{O_M^{\times }}{\nearrow } & \downarrow _{incl} &  & \quad
\quad \downarrow _{I_{*}\otimes Id} &  &  \\ 
H(K) &  & H((M,M\backslash K)\times K) & ^{\underrightarrow{\ \quad \psi
_{*}\quad \ }} & H(M^{\times })\otimes H(V) & ^{\underrightarrow{~m(P\otimes
Id)~}} & H(V), \\  
& \stackrel{O_K^{\times }}{\searrow } & \uparrow ^{incl} & \stackrel{}{%
\nearrow } &  &  &  \\  
&  & H((V,V\backslash K)\times K) &  &  &  & 
\end{array}
$$
where the vertical arrows are induced by inclusions and $\psi $ is given by%
$$
\psi (v,k)=(v,k,k),\qquad v\in M,\ k\in K, 
$$
so $\psi =(I\times Id)(Id\times \delta ).$ The diagram is commutative,
because the left triangle commutes by Definition \ref{def-fund}$.$ Now,
according to Lemma \ref{prelD}, $\psi _{*}:H((V,V\backslash K)\times
K)\rightarrow H(M^{\times }\times V)$ is also induced by%
$$
\psi ^{\prime }(v,k)=(v,k,v),\qquad v\in V,\ k\in K, 
$$
so $\psi ^{\prime }=(Id\times t)(\delta \times Id).$ Then the lower path
defines $\Phi .$ Therefore so does the upper one. Hence%
$$
\Phi =m(P\otimes Id)\psi _{*}O_M^{\times }. 
$$

Consider $u\in H_i(K)$ and $w=\delta _{*}(u)\in (H(K)\otimes H(V))_i$, then 
$$
w=\sum_{k+l=i}a_k\otimes b_l,\qquad a_k\in H_k(K),\ b_l\in H_l(V). 
$$
But $(\eta \otimes Id)\delta _{*}(u)=u,$ where $\eta :H(K)\longrightarrow 
{\bf Q}$ is the augmentation. Then $u=(\eta \otimes Id)\delta _{*}(u)=(\eta
\otimes Id)w=(\eta \otimes Id)\sum_{k+l=i}a_k\otimes b_l=\eta (a_0)\otimes
b_i.$ Therefore, $b_i=i_{*}(u)$ and $a_0$ is represented by some $p\in K.$
Since $O_M\in H_n(M)$, we have%
$$
\begin{array}[t]{ll}
\Phi (u)=m(P\otimes Id)\psi _{*}(O_M\otimes u) & =m(PI_{*}\otimes
Id)(O_M\otimes w) \\  
& =m(PI_{*}\otimes Id)(O_M\otimes \sum_{k+l=i}(a_k\otimes b_l)) \\  
& =m(\sum_{k+l=i}PI_{*}(O_M\otimes a_k)\otimes b_l) \\  
& =m(I_{*}(O_M\otimes a_0)\otimes b_i) \\  
& =I_{*}(O_M\otimes p)\cdot i_{*}(u).
\end{array}
$$
Finally, we observe that $l_pi_p(x)=I(x,p)$ for any $x\in M,\ \ p\in K,$ so%
$$
\begin{array}{ll}
\Phi (u) & =l_{p*}i_{p*}(O_M)\cdot i_{*}(u) \\  
& =l_{p*}(s(p))\cdot i_{*}(u)\quad 
\text{by Proposition \ref{ip(z)-sp}} \\  & =i_{*}(u)\quad \quad \quad \text{%
by Proposition \ref{lpsp}.}
\end{array}
$$
\end{proof}

In a fashion similar to the proof of Proposition \ref{I(fg)} the original
Dold's Lemma follows from this theorem.

\section{The Lefschetz Number of a Pair and Theorem \ref{identity}. \label
{MainRes}}

Now we recall some fundamentals of the theory of the Lefschetz number, see 
\cite[p. 207-208]{Dold} and \cite[p. 19-20]{Gorn}. Let%
$$
\begin{array}{l}
E_q^{*}=Hom(E_{-q}), 
\text{\quad }E^{*}=\{E_q^{*}\}, \\ (E^{*}\otimes E)_k=\bigotimes
{}_{q+i=k}(E_q^{*}\otimes E_i),\text{\quad }E^{*}\otimes E=\{(E^{*}\otimes
E)_k\}. 
\end{array}
$$
Now we define the following maps%
$$
\begin{array}{l}
e:(E^{*}\otimes E)_0\longrightarrow 
{\bf Q}\text{ by \quad }e(u\otimes v)=u(v)\text{ ({\it the evaluation map}),}
\\ \theta :(E^{*}\otimes E)_0\longrightarrow Hom(E,E)\text{ by \quad }[%
\theta (a\otimes b)](u)=(-1)^{\mid b\mid \cdot \mid u\mid }a(u)\cdot b, 
\end{array}
$$
where $|w|$ stands for the degree of $w.$

\begin{proposition}
\label{theta}\cite[Theorem II.1.5, p. 20]{Gorn} If $h:E\longrightarrow E$ is
an endomorphism of degree zero of a finitely generated graded module, then%
$$
e(\theta ^{-1}(h))=L(h). 
$$
\end{proposition}

What follows is an adaptation of the Gorniewicz's argument \cite[pp. 38-40]
{Gorn} to the new situation. The next two lemmas are trivial.

\begin{lemma}
\label{diag3} Let $J:H(V,V\backslash K)\longrightarrow (H(K))^{*}$ be a
homomorphism of degree $(-n)$ given by%
$$
J(u)(v)=I_{*n}(u\otimes v),\text{\quad }u\in H(V,V\backslash K),\ v\in H(K). 
$$
Then $I_{*n}=e(J\otimes Id),$ so that the following diagram commutes 
$$
\begin{array}{rrr}
(H(V,V\backslash K)\otimes H(K))_n & \stackrel{J\otimes Id}{\longrightarrow }
& ((H(K))^{*}\otimes H(K))_0 \\ 
\downarrow ^{I_{*}} & \swarrow _e &  \\ 
H_n(M^{\times })\simeq {\bf Q}. &  & 
\end{array}
$$
\end{lemma}

\begin{lemma}
\label{defa} Let 
$$
a=(J\otimes Id)(Id\otimes \varphi )\delta _{*}(O_K)=(J\otimes \varphi
)\delta _{*}(O_K)\in ((H(K)^{*}\otimes H(K))_0. 
$$
Then%
$$
e(a)=I_{*}(Id\otimes \varphi )\delta _{*}(O_K). 
$$
\end{lemma}

\begin{lemma}
\label{lastdiag} Let $\varphi :H(V)\rightarrow H(K)$ be a homomorphism. Then
the following diagram commutes:%
$$
\begin{array}{ccc}
H(V,V\backslash K)\otimes H(V)\otimes H(K) & ^{\underrightarrow{J\otimes
\varphi \otimes Id}} & (H(K))^{*}\otimes H(K)\otimes H(K) \\ 
~~~\downarrow ^{Id\otimes t_{*}} &  & ~~~~\downarrow ^{Id\otimes t_{*}} \\ 
H(V,V\backslash K)\otimes H(K)\otimes H(V) & ^{\underrightarrow{J\otimes
Id\otimes \varphi }} & (H(K))^{*}\otimes H(K)\otimes H(K) \\ 
~~~~~\downarrow ^{PI_{*}\otimes Id} &  & ~~~\downarrow ^{e\otimes Id} \\ 
H_n(M^{\times })\otimes H(V) &  & {\bf Q}\otimes H(K) \\ \downarrow ^m &  & 
\downarrow ^m \\ 
H(V) & ^{\underrightarrow{\qquad \varphi \qquad }} & H(K).
\end{array}
$$
\end{lemma}

\begin{proof} The first square trivially commutes. For the second, consider going $%
\rightarrow \downarrow $, then we get$:$%
$$
\begin{array}{cl}
m(e\otimes Id)(J\otimes Id\otimes \varphi ) & =m(e(J\otimes Id)\otimes
\varphi ) \\  
& =m(I_{*n}\otimes \varphi )\qquad 
\text{by Lemma }\ref{diag3} \\  & =I_{*n}\cdot \varphi \qquad \qquad \text{%
by definition of }m.
\end{array}
$$
Going $\downarrow \rightarrow ,$ we get%
$$
\begin{array}{ll}
\varphi m(PI_{*}\otimes Id) & =\varphi m(I_{*n}\otimes Id)\qquad 
\text{ by definition of }P \\  & =\varphi (I_{*n}\cdot Id)\qquad \qquad 
\text{by definition of }m \\  & =I_{*n}\cdot \varphi \qquad \qquad \qquad 
\text{by linearity of }\varphi .
\end{array}
$$
\end{proof}

The following theorem implies Theorem \ref{identity}.

\begin{theorem}
\label{L(fi)::fi(Ok)}Suppose $V$ is an ANR, $K$ is an arcwise connected
space, $H(K)$ is finitely generated. Then for any homomorphism $\varphi
:H(V)\rightarrow H(K)$ we have 
$$
L(\varphi i_{*})=I_{*}(Id\otimes \varphi )\delta _{*}(O_K), 
$$
where $i:K\longrightarrow V$ is the inclusion.
\end{theorem}

\begin{proof} We start in the left upper corner of the above diagram with $\delta
_{*}(O_K)\otimes u$, where $u\in H(K).$ Then going $\downarrow \rightarrow ,$
we get $\varphi i_{*}(u)$ by Theorem \ref{MDold}. Going $\rightarrow
\downarrow ,$ we get%
$$
\begin{array}{l}
\quad \ m(e\otimes Id)(Id\otimes t_{*})(J\otimes \varphi \otimes Id)(\delta
_{*}(O_K)\otimes u) \\ 
=m(e\otimes Id)(Id\otimes t_{*})((J\otimes \varphi )\delta _{*}(O_K)\otimes
u) \\ 
=m(e\otimes Id)(Id\otimes t_{*})(a\otimes u)\qquad 
\text{by Lemma }\ref{defa} \\ =m(e\otimes Id)(Id\otimes
t_{*})(\sum\nolimits_ia_i\otimes a_i^{\prime }\otimes u)\ 
\text{if }a=\sum\nolimits_ia_i\otimes a_i^{\prime },a_i\in
(H(K))^{*},a_i^{\prime }\in H(K) \\ =m(e\otimes
Id)\sum\nolimits_i(-1)^{\mid a_i^{\prime }\mid \cdot \mid u\mid
}(a_i\otimes u\otimes a_i^{\prime })\qquad  \\ 
=\sum\nolimits_im((-1)^{\mid a_i^{\prime }\mid \cdot \mid u\mid
}e(a_i\otimes u)\otimes a_i^{\prime }) \\ 
=\sum\nolimits_im((-1)^{\mid a_i^{\prime }\mid \cdot \mid u\mid
}a_i(u)\otimes a_i^{\prime })\qquad \qquad 
\text{by definition of }e \\ =\sum\nolimits_i(-1)^{\mid a_i^{\prime }\mid
\cdot \mid u\mid }a_i(u)\cdot a_i^{\prime }\qquad \qquad 
\text{by definition of }m \\ =\sum\nolimits_i\theta (a_i\otimes a_i^{\prime
})(u)\qquad \qquad 
\text{by definition of }\theta  \\ =\theta (a)(u).
\end{array}
$$
Thus $\theta (a)=\varphi i_{*}:H(K)\longrightarrow H(K).$ Since $H(K)$ is a
finitely generated graded module, Proposition \ref{theta} applies and we have%
$$
L(\varphi i_{*})=L(\theta (a))=e(a). 
$$
Now the statement follows from Lemma \ref{defa}. \end{proof}

\section{A Lefschetz-Type Coincidence Theorem for Maps to an Open Subset of
a Manifold.\label{GenCase2}}

Recall that if $(V,V\backslash K)$ is a manifold then according to Theorem 
\ref{representation} $\varphi =g_{*}f_{!}$ satisfies an identity that
connects it to the coincidence index. Our next goal is to show that even if $%
V$ is not a manifold under certain purely homological conditions we can
construct $\varphi $ satisfying that identity. This leads to the proof of a
Lefschetz-type theorem for Case 2.

\begin{proposition}
\label{degree} Suppose $f:(X,X\backslash N)\rightarrow (V,V\backslash K),$ $%
g:X\rightarrow K$ are continuous maps, and there is$\ \mu \in
H_n(X,X\backslash N)$ satisfying:

\begin{description}
\item[(a)] $f_{*}(\mu )=O_K,$
\end{description}

\noindent where $f_{*}:H_n(X,X\backslash N)\rightarrow H_n(V,V\backslash K)$%
, and there is a homomorphism $\varphi :H(V)\rightarrow H(K)$ of degree $0$
satisfying:

\begin{description}
\item[(b)] $\varphi f_{*}=g_{*},$
\end{description}

\noindent where $f_{*}:H(X)\rightarrow H(V)$.

\noindent Then the following holds%
$$
I_{fg}=I_{*}(Id\otimes \varphi )\delta _{*}(O_K), 
$$
i.e., $I_{fg}$ is the image of $O_K$ under the composition of the following
maps:%
$$
H(V,V\backslash K)^{\underrightarrow{\ \delta _{*}\ }}H(V,V\backslash
K)\otimes H(V)^{\underrightarrow{Id\otimes \varphi \ \ }}H(V,V\backslash
K)\otimes H(K)^{\underrightarrow{\ I_{*}\ }}H(M^{\times }). 
$$
\end{proposition}

\begin{proof} (cf. Gorniewicz \cite[pp. 15-16]{Gorn}) The following diagram
commutes:%
$$
\begin{array}{ccccc}
H(V,V\backslash K)\otimes H(V) & ^{\underleftarrow{Id\otimes f_{*}}} & 
H(V,V\backslash K)\otimes H(X) & ^{\underrightarrow{Id\otimes g_{*}}} & 
H(V,V\backslash K)\otimes H(K) \\ 
\uparrow ^{\alpha _1} &  & \uparrow ^{\alpha _2} &  & \uparrow ^{\alpha _3}
\\ 
H((V,V\backslash K)\times V) & ^{\underleftarrow{(Id\times f)_{*}}} & 
H((V,V\backslash K)\times X) & ^{\underrightarrow{(Id\times g)_{*}}} & 
H((V,V\backslash K)\times K) \\ 
\uparrow ^{\delta _{*}} &  & ~~~~\uparrow ^{(f,Id)_{*}} &  & \downarrow
^{I_{*}} \\ 
H(V,V\backslash K) & ^{\underleftarrow{\quad f_{*}\ \ }} & H(X,X\backslash
N) & ^{\underrightarrow{\ (f,g)_{*}}} & H(M^{\times }),
\end{array}
$$
where $\alpha _1,\ \alpha _2,\ \alpha _3$ are the isomorphisms from the
K\"unneth theorem. If we start with $\mu \in H(X,X\backslash N)$ in the
middle of the lower row of the diagram, then from commutativity of the
diagram it follows that 
$$
\begin{array}{ll}
I_{*}\alpha _3^{-1}(Id\otimes \varphi )\alpha _1\delta _{*}(O_K) & 
=I_{*}\alpha _3^{-1}(Id\otimes \varphi )\alpha _1\delta _{*}f_{*}(\mu )\quad 
\text{by (a)} \\  & =I_{*}\alpha _3^{-1}(Id\otimes \varphi )(Id\otimes
f_{*})\alpha _2(f,Id)_{*}(\mu )\quad 
\text{the left half} \\  & =I_{*}\alpha _3^{-1}(Id\otimes \varphi
f_{*})\alpha _2(f,Id)_{*}(\mu ) \\  
& =I_{*}\alpha _3^{-1}(Id\otimes g_{*})\alpha _2(f,Id)_{*}(\mu )\quad 
\text{by (b)} \\  & =I_{fg}\quad \text{the right half of the diagram. }
\end{array}
$$
\end{proof}

It is obvious that both (a) and (b) are satisfied when $f$ induces
isomorphisms $H(X)\simeq H(V)$ and $H_n(X,X\backslash N)\simeq
H_n(V,V\backslash K).$ Another example: if the maps $f,g:X\rightarrow M$
satisfy $f_{*n}\neq 0,\ g_{*}=0$, then we can select $\varphi =0$ so that $%
L(\varphi )=1$.

Now Proposition \ref{degree}\ and Theorem \ref{L(fi)::fi(Ok)} imply the
following.

\begin{theorem}[Lefschetz-Type Theorem]
\label{main} Suppose $X$ is a topological space, $N\subset X,$ $M$ is an
oriented connected compact closed $n$-manifold, $V\subset M$ is an ANR, $%
K\subset V$ is an arcwise connected space, $H(K)$ is finitely generated, $%
(M,V,M\backslash K)$ is an excisive triad. Suppose maps 
$$
f:(X,X\backslash N)\longrightarrow (V,V\backslash K),\text{\quad }%
g:X\longrightarrow K, 
$$
and a homomorphism $\varphi :H(V)\rightarrow H(K)$ satisfy the conditions of
Proposition \ref{degree}. Then the coincidence index is equal to the
Lefschetz number:%
$$
I_{fg}=L(\varphi i_{*}). 
$$
Moreover, if $L(\varphi i_{*})\neq 0,$ then $(f,g)$ has a coincidence.
\end{theorem}

Remark: We proved the theorem for singular homology, but the proof
is algebraic except for generalized Dold's Lemma \ref{MDold}. And since it
holds for \v Cech homology, then so does the theorem.

The following two examples show limits of applicability of this result.

\begin{example}
\label{Example0}(Dranishnikov \cite[Lemma 1.9]{Dran}) There is a multivalued
u.s.c. retraction $\Phi :{\bf D}^n\rightarrow {\bf S}^{n-1}:$%
$$
\Phi (x)=\{y\in {\bf S}^{n-1}:|y-x|\geq 4|x|^2-3|x|\}. 
$$
\end{example}

\begin{example}
\label{Example}Let ${\bf M}^2$ be the M\"obius band, given in cylindrical
coordinates by: $z=\theta ,\ -1\leq r\leq 1,\ 0\leq \theta \leq \pi ,$ with
the top and bottom edges identified. Let $f:({\bf M}^2,\partial {\bf M}%
^2)\rightarrow ({\bf D}^2,{\bf S}^1)$ be the projection on the horizontal
plane and $g:{\bf M}^2\rightarrow {\bf S}^1$ the projection on the $z$-axis.
\end{example}

Observe that $\Phi $ in Example \ref{Example0} has no fixed points, while in
Example \ref{Example} $g$ is homotopic to a map $g^{\prime }$ such that the
pair $(f,g^{\prime })$ has no coincidence. This is reflected in the fact
that $\Phi (x)$ fails to be acyclic for $|x|\leq 1/2,$ while in Example \ref
{Example} $f_{*}$ does not satisfy condition (a).

\section{The Generalized Lefschetz Number and Case 2.\label{VietorisMap}}

In Corollary \ref{GornCoin} below we will see how one can use the
Vietoris-Begle Theorem \ref{VietB} to avoid the restriction on relative
behavior of $f:(X,X\backslash N)\rightarrow (V,V\backslash K)$ and consider
only $f:X\rightarrow V.$ For this purpose, we would like to be able to deal
with the Lefschetz number of $i_{*}\varphi _{fg}:H(V)\longrightarrow H(V),$
instead of $\varphi _{fg}i_{*}:H(K)\longrightarrow H(K),$ as before. Then,
if $H(V)$ is not finitely generated, we need to define the generalized
Lefschetz number $\Lambda (\cdot )$, as in \cite[pp. 20-23]{Gorn}.

Let $h:E\rightarrow E$ be an endomorphism of an arbitrary vector space $E$.
Denote by $h^{(n)}:E\rightarrow E$ the $n$th iterate of $h$, then the kernels%
$$
\ker h\subset \ker h^{(2)}\subset \ldots \subset \ker h^{(n)}\subset \ldots 
$$
form an increasing sequence of subspaces of $E$. Let 
$$
N(h)=\bigcup\limits_n\ker h^{(n)}\text{ and }\tilde E=E/N(h). 
$$
Then $h$ induces the endomorphism $\tilde h:\tilde E\rightarrow \tilde E.$

\begin{definition}
Let $h=\{h_q\}$ be an endomorphism of degree $0$ of a graded vector space $%
E=\{E_q\}$ and suppose $\tilde E$ is finitely generated. Then the {\it %
generalized Lefschetz number (in the sense of Leray)} of $h$ is given by%
$$
\Lambda (h)=\sum\limits_q(-1)^qtr(\tilde h_q)=L(\tilde h). 
$$
\end{definition}

\begin{proposition}
\cite[II.2.3, p. 22]{Gorn} If $E$ is a finitely generated graded vector
space then $\Lambda (h)=L(h).$
\end{proposition}

\begin{proposition}
\label{genLef}\cite[II.2.4, p. 22]{Gorn} Let $E,E^{\prime }$ be graded
modules and suppose that the following diagram commutes:%
$$
\begin{array}{ccc}
E & ^{\underrightarrow{\quad k\quad }} & E^{\prime } \\ 
~~\uparrow ^h & \nwarrow ^l & ~\uparrow ^{h^{\prime }} \\ 
E & ^{\underrightarrow{\quad k\quad }} & E^{\prime }.
\end{array}
$$
Then, if $\Lambda (h)$ is defined then so is $\Lambda (h^{\prime })$ and $%
\Lambda (h)=\Lambda (h^{\prime })$.
\end{proposition}

\begin{theorem}
\label{mainG}Under conditions of Theorem \ref{main}, 
$$
I_{fg}=\Lambda (i_{*}\varphi _{fg}). 
$$
\end{theorem}

\begin{proof} It follows from Proposition \ref{genLef} and commutativity of the
diagram:%
$$
\begin{array}{ccc}
H(K) & ^{\underrightarrow{\quad i_{*}\quad }} & H(V) \\ 
~~~~~~\uparrow ^{\varphi _{fg}i_{*}} & \nwarrow ^{\varphi _{fg}} & 
~~~~\uparrow ^{i_{*}\varphi _{fg}} \\ 
H(K) & ^{\underrightarrow{\quad i_{*}\quad }} & H(V).
\end{array}
$$
\end{proof}

This theorem implies the Coincidence Theorem of Gorniewicz \cite[p. 38]{Gorn}%
, as follows.

\begin{corollary}
\label{GornCoin}Suppose $X$ is a topological space, $V\subset {\bf R}^n$ is
open. Suppose%
$$
f,g:X\rightarrow V 
$$
are two continuous maps such that $f$ is Vietoris and $g$ is compact (i.e., $%
\overline{g(X)}$ is compact). Then if $\Lambda (g_{*}f_{*}^{-1})\neq 0$ with
respect to \v Cech homology over ${\bf Q}$, then the pair $(f,g)$ has a
coincidence.
\end{corollary}

\begin{proof} Since $g$ is a compact map, there is a finite connected polyhedron $K$
such that $g(X)\subset K\subset V$. We can assume that $V\subset M={\bf S}^n$%
. Then $(M,V,V\backslash K)$ is an excisive triad, $V$ is an ANR, $K$ is
arcwise connected space. Let $N=f^{-1}(K)$, then $Coin(f,g)\subset N$. Since 
$f_{*}$ is an isomorphism by Proposition \ref{VietB}, all the conditions of
Proposition \ref{degree} are satisfied. Therefore by the theorem, $%
I_{fg}=\Lambda (g_{*}f_{*}^{-1}).$ \end{proof}

See Gorniewicz \cite[pp. 40-43]{Gorn} for applications of this theorem to
the study of fixed points of multivalued maps on polyhedra, ANRs, etc.\\

{\bf Acknowledgments.} This paper is a version of Chapter I of my doctoral thesis written under the guidance of M.-E. Hamstrom at the University of Illinois at Urbana-Champaign. I would like to express my deep gratitude to
Professor Robert F. Brown who read an early version of the manuscript,
suggested ways to turn it into a publication, and always provided advice and
encouragement. I also thank Professor R. McCarthy for valuable conversations
regarding this work and the referee for a very thorough and thoughtful
report.

\end{document}